\documentclass{amsart}

\usepackage{amsfonts}
\usepackage{amsmath}
\usepackage{amssymb}
\usepackage{latexsym}
\usepackage{graphicx}

\newtheorem{thm}{\bf Theorem}[section]
\newtheorem{cor}[thm]{\bf Corollary}
\newtheorem{lem}[thm]{\bf Lemma}
\newtheorem{prop}[thm]{\bf Proposition}
\newtheorem{defn}[thm]{\bf Definition}

\newtheorem{exmp}[thm]{\bf Example}

\catcode`\ç=13
\defç{\c{c}}
\catcode`\é=13
\defé{\'e}
\catcode`\à=13
\defà{\`a}
\catcode`\è=13
\defè{\`e}
\catcode`\â=13
\defâ{\^a}
\catcode`\ù=13
\defù{\`u}
\catcode`\ê=13
\defê{\^e}
\catcode`\î=13
\defî{\^\i}
\catcode`\ô=13
\defô{\^o}
\newcommand{\field}[1]{\mathbb{#1}}

\newcommand{\N }{\field{N}}

\def\proof{{\parindent0pt {\bf Proof.\ }}}

\def\wdim{{\rm wdim}}
\def\gldim{{\rm gldim}}

\def\Gwdim{{\rm G\!-\!wdim}}
\def\Ggldim{{\rm G\!-\!gldim}}

\def\GPD{{\rm GPD}}
\def\GID{{\rm GID}}

\def\pd{{\rm pd}}
\def\fd{{\rm fd}}
\def\id{{\rm id}}

\def\Gpd{{\rm Gpd}}
\def\Gfd{{\rm Gfd}}
\def\Gid{{\rm Gid}}

\def\Im{{\rm Im}}

\def\Ext{{\rm Ext}}

\def\Hom{{\rm Hom}}

\def\sup{{\rm sup}}
\def\max{{\rm max}}

\newcommand{\cqfd}
{\hspace{1cm}
\rule{2mm}{2mm}%
\medbreak%
\par%
}
\def\1{{\noindent\rm (1)}}
\def\2{{\noindent\rm (2)}}
\def\3{{\noindent\rm (3)}}
\def\4{{\noindent\rm (4)}}
\def\5{{\noindent\rm (5)}}
\begin{document}

\title[Global Gorenstein dimensions of polynomial rings and
 ...]{Global Gorenstein dimensions of
polynomial rings and of direct product of rings}

\author{D. Bennis}
\address{Department of Mathematics, Faculty of Sciences and Technology, University S. M. Ben Abdellah, Fez 30000, Morocco}
\email{driss\_bennis@hotmail.com}

\author{N. Mahdou}
\address{Department of Mathematics,Faculty of Sciences and Technology, University S. M. Ben Abdellah, Fez 30000, Morocco}
\email{mahdou@hotmail.com}


\begin{abstract}In this paper, we extend the well-known  Hilbert's
syzygy theorem to the Gorenstein homological dimensions of rings.
Also, we study the Gorenstein homological dimensions of direct
product of rings, which gives examples of non-Noetherian rings of
finite Gorenstein dimensions and infinite classical weak
dimension.
\end{abstract}
\maketitle
\begin{section}{Introduction}  Throughout this paper all
rings are commutative with identity element and all modules are
unital.\bigskip

\textbf{Setup and Notation:} Let  $R$ be a ring, and let $M$ be an
$R$-module.\\
As usual we use $\pd_R(M),\ \id_R(M)$, and $\fd_R(M)$
  to denote, respectively, the classical projective, injective and flat
dimensions of $M$. $\gldim(R)$ and $\wdim(R)$ are, respectively,
the classical global and weak dimensions of $R$.\\
It is by now a well-established fact that even if $R$ to be
non-Noetherian, there exist  Gorenstein projective, injective and
flat dimensions of $M$, which are   usually denoted by
$\Gpd_R(M),\ \Gid_R(M)$, and $\Gfd_R(M)$, respectively. Some
references are  \cite{BM, BM2,LW, CFH, HH}.\bigskip

Recently, the authors in \cite{BM2} started the study of global 
Gorenstein dimensions of rings, which are called, for
a ring $R$, Gorenstein projective, injective, and weak dimensions
of $R$, denoted  by  $\GPD(R)$,  $\GID(R)$, and $G-wdim(R)$,
respectively, and, respectively, defined  as follows:\\

  \label{D-Gdim-ring}

 \noindent   \mbox{}\hspace{3cm}  $\GPD(R) = \sup\{\Gpd_R(M)\,|\,M\;R\!-\!module\}$,\\
             \mbox{}\hspace{3cm}   $ \GID(R) = \sup\{\Gid_R(M)\,|\,M\;R\!-\!module\}$,    and\\
             \mbox{}\hspace{3cm} $G-wdim(R) = \sup\{\Gfd_R(M)\,|\,M\;R\!-\!module\}$.\\

They proved that, for any ring  $R$: $ G-wdim(R)\leq \GID(R)=
\GPD(R)$ \cite[Theorems 3.2 and 4.2]{BM2}. So, according to
the terminology of the classical theory of homological dimensions
of rings,  the common value of $\GPD(R)$ and $\GID(R)$ is called
\textit{Gorenstein global dimension} of $R$, and denoted by 
\boldmath $\mathbf{\Ggldim(R)}$. \unboldmath \\ 

 They also proved that the Gorenstein  global and weak dimensions
  are refinements of the classical   global and weak  dimensions of rings, respectively.
  That is \cite[Propositions 3.11 and 4.5]{BM2}:
  $\Ggldim(R)\leq \gldim(R)$ and
$\Gwdim(R)\leq \wdim(R)$, with equality if $\wdim(R)$ is
finite.\bigskip

In Section 2, we are mainly interested in computing  the
Gorenstein dimensions  of commutative polynomial rings to
generalize the fact that a polynomial ring over Gorenstein ring is
also Gorenstein. Explicitly,  we extent the equalities of the
well-known Hilbert's syzygy theorem to the Gorenstein global and
weak dimensions (please see Theorems \ref{GPsyzygy} and
\ref{GFsyzygy}).\bigskip

Recall that a ring $R$ is said to be
 $n$-Gorenstein, for a positive integer $n$, if it is Noetherian
 and $\id_R(R)\leq n$ \cite{Iwa}. And $R$ is said to be Gorenstein, if
it is $n$-Gorenstein for some positive integer $n$. These kind of
rings have a characterization by the weak and global Gorenstein 
dimensions of rings. Namely, we have \cite[Corollary 2.3]{BM2}: If
$R$ is a  Noetherian ring, then $\Gwdim(R)=\Ggldim(R)$, such that,
for a positive integer $n$, $\Ggldim(R)\leq n$ if, and only if, $
R$ is  $n$-Gorenstein. In the case $n=0$, we have the last
equivalence  without assuming in the first place that $R$ to be
Noetherian; that is \cite[Theorem 5.2.4]{BM2}: $\Ggldim(R)= 0 $
if, and only if, $R$ is  $0$-Gorenstein\ (i.e.,\ quasi-Frobenius).
Thus, the study of the Gorenstein dimensions of polynomial rings
allows us to give examples of Noetherian rings which have finite
Gorenstein global (=weak) dimension and infinite
 global (=weak) dimension (see
Example \ref{Exm1}). Obviously, one would like to have examples
out of the class of Noetherian rings. Explicitly, one would like
to have examples  of a family of non-Noetherian rings
$\{R_{i}\}_{i\geq 1}$ such that $\Gwdim(R_i)=i$, $\Ggldim(R_i)>i$,
and $\wdim(R_i)=\infty$ for all
 $i\geq 1$. This will be obtained in Section 3 (Example
\ref{Exm2}) after a study of the Gorenstein dimensions of direct
product of rings (Theorems \ref{GP-direct-product} and
\ref{GF-direct-product}).
\end{section}
\begin{section}{Gorenstein dimensions of polynomial rings}
The aim of this section is to compute the Gorenstein
 global and weak  dimensions of commutative polynomial rings.\bigskip

Through this section,  $R$ denotes  a commutative ring, and $R[X]$
denotes the polynomial ring in one indeterminate over  $R$. We use
$M[X]$, where $M$ is an $R$-module, to denote the $R[X]$-module
$M\otimes _R R[X]$.

  The first main result in this
section is:
\begin{thm}\label{GPsyzygy}
Let  $R[X_1,X_2,...,X_n]$ be the polynomial ring in $n$
indeterminates over    $R$. Then:
 \boldmath $\mathbf{\Ggldim(R[X_1,X_2,...,X_n])=\Ggldim(R)+n}.$\unboldmath
\end{thm}
To prove this theorem, we need some results.\bigskip

First, we need the notion of strongly Gorenstein projective modules,
which are introduced in \cite{BM} to characterize the Gorenstein
projective modules.

\begin{defn}[\cite{BM}, Definition 2.1]\label{def-G-str-proj-mod}
A module $M$ is said to be strongly Gorenstein projective, if
there exists an exact sequence of the form:
 $ \mathbf{P}=\
\cdots\stackrel{f}{\longrightarrow}P\stackrel{f}{\longrightarrow}P\stackrel{f}{\longrightarrow}P
\stackrel{f}{\longrightarrow}\cdots$
     where $P$ is a projective $R$-module and $f$ is an endomorphism of $P$,
     such that  $M \cong \Im(f)$ and such that $\Hom ( -, Q) $ leaves the
sequence $\mathbf{P}$ exact whenever $Q$ is a projective module.
\end{defn}

These strongly Gorenstein projective  modules are a simple
characterization, and they are used to characterize the Gorenstein
projective modules. We have:

\begin{prop}[\cite{BM2}, Proposition 3.9]\label{Cara-G-str-proj}  A module $M$ is  strongly Gorenstein projective if, and only if,
there exists a short exact sequence of modules $0\rightarrow
M\rightarrow P\rightarrow M\rightarrow 0$ where $P$ is a
projective  module, and  $\Ext^i(M,Q)=0$ for some integer  $i>0$
and for any  module $Q$  with finite projective dimension (or for
any projective module $Q$).
\end{prop}

\begin{thm}[\cite{BM}, Theorem 2.7] \label{tcara-G-stro-pro}
A module is Gorenstein projective  if, and only if, it is a direct
summand of a strongly Gorenstein projective module.
\end{thm}

\begin{lem}\label{GMX-leq-GM}
For any  $R$-module $M$,
 $\Gpd_{R[X]}(M[X])\leq \Gpd_{R}(M).$
\end{lem}
 \proof The standard argument shows that the inequality is a simple
consequence of the implication: $M$ is a Gorenstein projective
$R$-module $\Rightarrow$ $M[X]$ is a Gorenstein projective
$R[X]$-module.
 So, assume  $M$ to be a Gorenstein projective $R$-module. We claim
 that $M[X]$ is a Gorenstein projective
$R[X]$-module.\\
 From  Theorem \ref{tcara-G-stro-pro}, $M$ is a direct summand of a strongly
 Gorenstein projective $R$-module  $N$. Then, by
 Proposition  \ref{Cara-G-str-proj}, there exists a short exact
sequence of $R$-modules
 $0\rightarrow N\rightarrow P\rightarrow N\rightarrow 0$
where $P$ is a projective $R$-module  and $\Ext_{R}(N,Q)=0$
 for any projective $R$-module  $Q$.
Hence, we have the following  short exact sequence of
$R[X]$-modules
 $0\rightarrow N[X]\rightarrow P[X]\rightarrow N[X]\rightarrow 0$
such that  $P[X]$ is a projective $R[X]$-module; and from
\cite[Theorem 11.65]{Rot}, $\Ext_{R[X]}(N[X],Q)\cong
\Ext_{R}(N,Q)=0$ for every projective $R[X]$-module $Q$ (since $Q$
is also projective $R$-module). Thus, $N[X]$ is a  strongly
      Gorenstein projective    $R[X]$-module (by Proposition  \ref{Cara-G-str-proj}). Therefore,  $M[X]$ is a Gorenstein projective $R[X]$-module as a
        direct summand of the strongly Gorenstein projective $R[X]$-module
        $N[X]$ (by  \cite[Theorem 2.5]{HH}).\cqfd

\begin{lem}\label{ineGpddd}
 For any short exact sequence  of modules $0\rightarrow
A\rightarrow B\rightarrow C\rightarrow 0$, we have inequality:
 $\Gpd(C)\leq 1+\max\{\Gpd(A),\Gpd(B)\}.$ \\

\end{lem}
\proof  The argument is analogous   to the proof of
\cite[Corollary 1.2.9 (a)]{LW}.\cqfd

\begin{lem}\label{ineGpddd}
 Let $R[X]$ be the polynomial ring in one indeterminate over $R$. Then, 
  $G-gldim(R[X]) \leq  G-gldim(R)+1.$ \\

\end{lem}
\proof  Using Lemmas 2.6 and 2.5, the argument is analogous  to the proof from 
\cite[Lemma 9.29]{Rot}.\cqfd

\begin{lem}\label{GMX=GM} If all projective  $R[X]$-modules  have
finite injective dimension, then, for any $R$-module $M$,
$Gpd_{R[X]}(M[X])= \Gpd_{R}(M).$\\
In particular, $G-gldim(R) \leq  G-gldim(R[X]).$ \\
\end{lem}

\proof From Lemma \ref{GMX-leq-GM}, it remains to prove: $\Gpd_{R[X]}(M[X]) \geq \Gpd_{R}(M).$\\
The standard argument shows that this inequality   follows from
the implication:  $M[X]$ is a Gorenstein projective $R[X]$-module
$\Rightarrow$ $M$  is a Gorenstein projective $R$-module. So,
assume that $M[X]$ is  a  Gorenstein projective $R[X]$-module. We
claim
 that $M$ is a Gorenstein projective
$R$-module.\\ From Theorem \ref{tcara-G-stro-pro},  $M[X]$   is a
direct summand of a strongly Gorenstein projective $R[X]$-module
$N$. For such $R[X]$-module  $N$, and from Proposition
  \ref{Cara-G-str-proj}, there exists a short exact sequence of
$R[X]$-modules
 $(\star)\  0\rightarrow N\rightarrow P\rightarrow
N\rightarrow 0$  where $P$ is a projective $R[X]$-module (then it
is also projective as an $R$-module), and
$\Ext_{R[X]}(N,Q')=0$ for every  $R[X]$-module $Q'$ with finite projective dimension.\\
Note that $M$ is a direct summand of $N$, since it is a direct
summand of $M[X]$. Then, to conclude the proof, it suffices, by
\cite[Theorem 2.5]{HH}, to prove that  $N$ is also a strongly
Gorenstein projective $R$-module. For that, and from the short
exact sequence $(\star)$, it remains to prove, from Proposition
\ref{Cara-G-str-proj}, that $\Ext^{n}_{R}(N,L)=0$
  for some integer $n>0$ and for every  projective $R$-module $
  L$. From \cite[Theorem 11.65]{Rot}, $\Ext^{i}_{R[X]}(N[X],L)\cong
\Ext^{i}_{R}(N,L)$  for every $R[X]$-module $L$. If $L$ is a
projective $R$-module, then $\pd_{R[X]}(L)=1$ (by  \cite[Theorem
C, page 124]{Kap}). Thus, the hypothesis  implies that
$\id_{R[X]}(L)$ is finite. Then, there exists an integer $n>0$
such that $\Ext^{n}_{R[X]}(N[X],L)=0$. Therefore,
$\Ext^{n}_{R}(N,L)=0$, as desired. \\
Finally, $G-gldim(R) \leq  G-gldim(R[X])$ is clear by above and this completes the 
proof of Lamma 2.8. \cqfd

 \bigskip

\noindent\textbf{Proof of Theorem \ref{GPsyzygy}.} By induction on
$n$, we may assume $n=1$, such that we write $R[X]=R[X_1]$.
Lemmas 2.7 and 2.8 give: $$ \Ggldim(R)\leq
\Ggldim(R[X])\leq \Ggldim(R)+1. $$ Then, we may assume that
$\Ggldim(R)$ and $\Ggldim(R[X])$ are finite.\\   Let
$\Ggldim(R)=m$, for some positive integer $m$. From \cite[Lemma
3.3]{BM2}, there exist an $R$-module $M$ and a projective
$R$-module $P$ such that: $\Ext^{m}_{R}(M,P) \neq 0 $. Then, from
Rees's theorem \cite[Theorem 9.37]{Rot}:
 $\Ext^{m+1}_{R[X]}(M,P[X])\cong\Ext^{m}_{R}(M,P)\neq 0$, which implies, from \cite[Theorem 2.20]{HH} and since $P[X]$ is a
projective $R[X]$-module, that $\Gpd_{R[X]}(M)\geq m+1$.
Therefore,
 $\Ggldim(R[X])\geq m+1$, and then the desired
equality.\cqfd

    \begin{cor}
Let $R[X_1,X_2,...,X_n,...]$ be the polynomial in  infinity of
indeterminates. Then, $\Ggldim(R[X_1,X_2,...,X_n,...])=\infty$.
\end{cor}

\proof Obvious.\cqfd\bigskip

 Theorem  \ref{GPsyzygy} allows us to give a family
$\{R_i\}_{i\in \N}$ of commutative rings, such that
$\Ggldim(R_i)=i$ and $\wdim(R_i)=\infty$ for all $i\in \N$, as
shown by the following example:

\begin{exmp}\label{Exm1}
Consider a non-semisimple quasi-Frobenius ring (for example, 
$K[X]/(X^{2})$, where $K$ is a field). Then, for every
positive
 integer $n$, we have:\\
 \mbox{}\hspace{0.5cm}  $\Ggldim(R[X_1,X_2,...,X_n])=n\ $ and $\ \wdim(R[X_1,X_2,...,X_n])=\infty$.
\end{exmp}
\proof  We have $\Ggldim(R)=0$ since $R$ is quasi-Frobenius. On the other hand, 
$gldim(R) (=wdim(R)) =\infty$ by \cite[Proposition 3.11]{BM2} since 
$gldim(R) \not= 0$ (since $R$ is a non-semisimple ring) and 
$R$ is Noetherian (since $R$ is quasi-Frobenius). \cqfd

\bigskip

Now, we treat the Gorenstein weak dimension of polynomial rings.

\begin{thm}\label{GFsyzygy}
If the polynomial ring $R[X]$ in one indeterminate $X$ over $R$
is coherent, then:
\boldmath$\mathbf{\Gwdim(R[X])=\Gwdim(R)+1}.$\unboldmath
\end{thm}

\proof  First, note that $R=R[X]/(X)$ is also coherent (by
\cite[Theorem 4.1.1 (1)]{Glaz}), since $R$ is a finitely presented
$R[X]$-module (by the short exact sequence of $R[X]$-modules
$0\rightarrow R[X]\rightarrow R[X] \rightarrow R\rightarrow0$),
and since $R[X]$ is   coherent.\\
We proceed similarly to the proof of Theorem \ref{GPsyzygy}. Then,
we first show that  $\Gwdim(R[X])\leq \Gwdim(R)+1.$ This
inequality is obtained in the same way as the one of Corollary
\ref{co3} using \cite[Proposition 3.10]{HH}  instead of Lemma
\ref{GMX-leq-GM} and the fact that the inequality of Lemma
\ref{ineGpddd} remains true, over coherent rings, for the
Gorenstein flat dimension case.\\
Secondly, we must show that $\Gwdim(R)\leq \Gwdim(R[X])$. For
that, it suffices to prove the inequality $\Gfd_{R[X]}(M[X]) \geq
\Gfd_{R}(M)$, for every $R$-module $M$. Also, the proof of this
inequality is similar to the one in the proof of  Lemma
\ref{GMX=GM}, using the properties of strongly Gorenstein flat
modules  \cite[Definition 3.1]{BM} (\cite[Proposition 4.6]{BM2}
and \cite[Theorem 3.5]{BM}) instead of the one of the strongly
Gorenstein projective modules used in the proof of Lemma
\ref{GMX=GM}, and using the fact that, over coherent rings, the
set of all Gorenstein flat modules is closed under direct summands
\cite[Theorem 3.7]{HH}. And finally, we use   \cite[Theorem
11.64]{Rot}, \cite[Theorem 202]{Kap}), and \cite[Theorem  4.11
(8)]{BM2} instead of \cite[Theorem 11.65]{Rot},
 \cite[Theorem C, page 124]{Kap}, and  \cite[Proposition
 3.13]{BM2}, respectively.\\
Now, we have  the inequality $\Gwdim(R)\leq \Gwdim(R[X]) \leq
 \Gwdim(R)+1$, then we may assume that
$\Gwdim(R)$ and $\Gwdim(R[X])$  are finite.\\
Let $\Gwdim(R)=m$, for some positive integer $m$.  Hence, from
\cite[Theorem 4.11]{BM2}, there exists a  finitely presented
$R$-module $M$ such that $\Ext^{m}_{R}(M,R)\neq 0$. Then, from
Rees's theorem \cite[Theorem 9.37]{Rot}:
 $\Ext^{m+1}_{R[X]}(M,R[X])\cong\Ext^{m}_{R}(M,R)\neq 0.$ The
 $R$-module
 $M$ is also finitely presented as an $R[X]$-module (by
\cite[Theorem 2.1.7]{Glaz}, and since $R$ is a finitely presented
$R[X]$-module, by the short exact sequence of $R[X]$-modules
$0\rightarrow R[X]\rightarrow R[X] \rightarrow R\rightarrow0$).
Then, from \cite[Theorem 4.11]{BM2},  $\Gwdim(R[X])\geq m+1$,
which implies the desired equality.\cqfd

\end{section}
\begin{section}{Gorenstein  dimensions of  direct
product of rings}

The aim of this section is to compute the Gorenstein  global and
weak  dimensions of  direct product of commutative rings.\bigskip

We begin by the Gorenstein global dimension case, which is a
generalization of the classical equality: $\gldim(\Pi_{i=1}^m
R_i)=\sup\{\gldim(R_i), 1\leq i \leq m\}$ where
$\{R_{i}\}_{i=1,...,m}$  is a family of rings \cite[Chapter VI,
Exercice 8, page 123]{H.Cartan}.

\begin{thm}\label{GP-direct-product}
Let $\{R_{i}\}_{i=1,...,m}$  be a family of rings. Then:
 \boldmath $$\mathbf{\Ggldim(\displaystyle\prod_{i=1}^m R_i)=\sup\{\Ggldim(R_i), 1\leq i \leq m\}}.$$\unboldmath
\end{thm}

To prove this theorem, we need the following results:

\begin{lem}\label{Ch-ring-SGp}
Let $R\rightarrow S$ be a ring homomorphism such that $S$ is a
projective $R$-module.  If $M$ is a (strongly) Gorenstein
projective $R$-module, then $M \otimes_R S$ is a (strongly)
Gorenstein  projective $S$-module.\\
Namely, we have:
  $\Gpd_{S}(M\otimes_R S)\leq \Gpd_{R}(M).$
\end{lem}
\proof Assume at first $M$ to be a strongly Gorenstein projective
$R$-module. Then, there exists a short exact sequence of
$R$-modules $0\rightarrow M \rightarrow P \rightarrow M
\rightarrow 0$ where $P$ is a projective $R$-module, and
$\Ext_R(M, Q) = 0$ for any projective $R$-module $Q$.\\
Since $S$ is a projective (then flat) $R$-module, we have a short
exact sequence of $S$-modules $0\rightarrow M \otimes_R
S\rightarrow P\otimes_R S \rightarrow M\otimes_R S \rightarrow 0$
such that  $ P\otimes_R S $ is a projective $S$-module, and for
any
  projective $S$-module (then  projective
$R$-module) $L$,
 $\Ext_S(M\otimes_R S , L)=\Ext_R(M,L)=0$ (by \cite[Theorem
 11.65]{Rot}).
This implies that $M\otimes_R S $ is a strongly Gorenstein projective $S$-module.\\
Now, let $M$ be any arbitrary Gorenstein projective $R$-module.
Then, it is a direct summand of a strongly Gorenstein projective
$R$-module $N$. Then, $M\otimes_R S $  is a direct summand of the
$S$-module  $N \otimes_R S$ which is, from the reason above,
 strongly Gorenstein projective. Therefore,    $M\otimes_R S $
 is a  Gorenstein projective $S$-module, as desired.\cqfd

\begin{lem}\label{lem-produ-SGp}
Let $\{R_{i}\}_{i=1,...,m}$  be a family of rings such that all
projective  $R_i$-modules  have finite injective dimension, for
i=1,...,m. Let $M_{i}$ be an $R_i$-module for i=1,...,m. If each
$M_i$ is a (strongly) Gorenstein projective $R_i$-module, then
$\prod_{i=1}^m M_i$ is a (strongly) Gorenstein projective
$(\prod_{i=1}^m
R_i)$-module.\\
Namely, we have:
 $\Gpd_{(\Pi_{i=1}^{m} R_i)}(\prod_{i=1}^{m} M_i)\leq \sup\{\Gpd_{ R_i}(M_i), 1\leq i \leq
m\}.$
\end{lem}
 \proof  By induction on $m$, it suffices to prove the assertion for
 $m=2.$\\
We assume at first that $M_i$ is a strongly Gorenstein projective
$R_i$-module for $i=1,\   2$. Then, there exists a  short exact
sequence of $R_i$-modules $0\rightarrow M_i \rightarrow P_i
\rightarrow M_i \rightarrow 0$ where $P_i$ is a projective
$R_i$-module. Hence, we have a   short exact sequence of
$R_1\times R_2$-modules $0\rightarrow M_1\times M_2\rightarrow
P_1\times P_2  \rightarrow M_1\times M_2\rightarrow 0$ where
$P_1\times P_2 $ is a projective $R_1\times R_2$-module (by  \cite[Lemma 2.5 (2)]{Mah2001}).\\

On the other hand, let $Q$ be a projective $R_1 \times
R_2$-module. We have: $$Q= Q \otimes_{R_1 \times R_2} (R_1 \times
R_2)= Q\otimes_{R_1 \times R_2} (R_1 \times0 \oplus 0 \times R_2)=
Q_1 \times Q_2$$ where $Q_i= Q \otimes_{R_1 \times R_2} R_i$
 for $i=1,\    2.$  From \cite[Lemma 2.5
(2)]{Mah2001}, $Q_i$ is a projective $R_i$-module for $i=1,\ 2$.
Hence, by hypothesis $\id_{R_i}(Q_i)< \infty$ for $i=1,\   2$, and
from \cite[Chapter VI, Exercice 10, page 123]{H.Cartan},
$\id_{R_1\times R_2}(Q_i)\leq \id_{R_i}(Q_i)< \infty$ for $i=1,\
2$. Thus, $\id_{R_1\times R_2}(Q_1\times Q_2)< \infty$, so
$\Ext_{R_1\times R_2}^{k}(M_1\times M_2,Q_1\times Q_2) = 0$ for
some positive integer $k$. This implies, from   Proposition
\ref{Cara-G-str-proj}, that $M_1\times M_2$ is a
strongly Gorenstein projective $R_1\times R_2$-module.\\
Now, let $M_i$ be any arbitrary Gorenstein projective $R_i$-module
for $i=1,\   2$. Then, there exist  an $R_i$-module $G_i$ and a
strongly Gorenstein projective  $R_i$-module $N_i$ for $i=1,\   2$
such that $M_i \oplus G_i = N_i$. Then, $(M_1 \times M_2)\oplus
(G_1 \times G_2)=(M_1\oplus G_1)\times (M_2\oplus G_2) = N_1\times
N_2$. Since, by the reason above, $N_1\times N_2$ is a strongly
Gorenstein projective $R_1\times R_2$-module, and from Theorem
\ref{tcara-G-stro-pro}, $M_1 \times M_2$ is a   Gorenstein
projective $R_1\times R_2$-module, as desired.\cqfd\bigskip

\noindent\textbf{Proof of Theorem \ref{GP-direct-product}.} By
induction on $m$, it suffices to prove the equality  for
 $m=2.$
 We must show that, for any positive integer $d$, we have the equivalence:
 $$\Ggldim(R_1 \times R_2) \leq d \Longleftrightarrow \Ggldim(R_1) \leq d
 \quad and\quad \Ggldim(R_2) \leq d.$$
  \noindent  Then, assume that $\Ggldim(R_1 \times R_2)
\leq d$ for some positive integer $d$.\\
Let $M_i$ be an $R_i$-module for $i=1,\   2$. Since each $R_i$ is
a projective $R_1 \times R_2$-module, and from Lemma
\ref{Ch-ring-SGp}, we have:
 $\Gpd_{R_i}(M_i)=\Gpd_{R_i}((M_1 \times M_2) \otimes_{R_1\times
R_2} R_i)\leq \Gpd_{R_1\times R_2}(M_1 \times M_2)\leq d.$ This
follows that $\Ggldim(R_i) \leq d$ for $i=1,\   2$.\\
Conversely, assume that $\Ggldim(R_i) \leq d$ for $i=1,\   2$
(where $d$ is a positive integer), and consider an $R_1 \times
R_2$-module $M$. We may write  $M=  M_1 \times M_2$  where $M_i= M
\otimes_{R_1 \times R_2} R_i$
 for $i=1,\   2.$  By hypothesis and from  \cite[Proposition 3.14]{BM2}, we may apply
Lemma \ref{lem-produ-SGp}, and so
  $\Gpd_{R_1 \times R_2} (M_1 \times M_2)\leq \sup\{\Gpd_{R_1}(M_1), \Gpd_{R_2}(M_2)\}\leq
  d$.   Therefore,  $\Ggldim(R_1 \times R_2) \leq d $.\cqfd\bigskip

We can now construct a family of non-Noetherian rings
$\{R_{i}\}_{i\geq 1}$ such that  $\Ggldim(R_i)=i$ and
$\wdim(R_i)=\infty$ for all $i\geq 1$, as follows:

\begin{exmp}\label{Exm1-non-Noeth}
Consider  a  non-semisimple quasi-Frobenius  ring $R$, and a
non-Noetherian hereditary ring $S$. Then, for
every positive integer $n$, we have:\\
 $\Ggldim(R\times S[X_1,X_2,...,X_n])=n+1$  and  $
   \wdim(R\times S[X_1,X_2,...,X_n])=\infty.$
\end{exmp}
\proof Since $S$ is  a non-Noetherian hereditary ring,
$\gldim(S)=1$. Then, the first equality follows immediately from
Hilbert's syzygy theorem, Theorem \ref{GP-direct-product}, and
since  $\Ggldim(R)=0$ (since $R$ is quasi-Frobenius).\\
Now, if $\wdim(R\times S[X_1,X_2,...,X_n])<\infty$, then, from
\cite[Proposition 3.11]{BM2}, we have:
 $\gldim(R\times S[X_1,X_2,...,X_n])=\Ggldim(R\times
S[X_1,X_2,...,X_n])=n+1.$ Thus,   $\gldim(R)\leq n+1$ (by
\cite[Chapter VI, Exercice 10, page 123]{H.Cartan}). But, this is
absurd, since $R$ is a non-semisimple quasi-Frobenius
ring.\cqfd\bigskip


Now we study the Gorenstein weak dimension case.

\begin{thm}\label{GF-direct-product}
Let $\{R_{i}\}_{i=1,...,m}$  be a family of coherent rings. Then:
 \boldmath$$\mathbf{\Gwdim(\displaystyle\prod_{i=1}^m R_i)=\sup\{\Gwdim(R_i), 1\leq i \leq m\}}.$$\unboldmath
\end{thm}

Similarly to the Gorenstein global dimension case, we need
  the following lemmas:

\begin{lem}\label{lem-produ-SGf}
Let $\{R_{i}\}_{i=1,...,m}$  be a family of coherent rings  such
that all injective  $R_i[X]$-modules  have finite flat dimension,
for i=1,...,m. Let $M_{i}$ be a family of $R_i$-module for
i=1,...,m. If each $M_i$ is a (strongly) flat $R_i$-module, then
$\prod_{i=1}^m M_i$ is a (strongly) flat $(\prod_{i=1}^m
R_i)$-module.\\
Namely, we have: $\Gfd_{(\Pi_{i=1}^{m} R_i)}(\prod_{i=1}^m
M_i)\leq \sup\{\Gfd_{ R_i}(M_i), 1\leq i \leq m\}.$
\end{lem}
\proof First, note  that $\Pi_{i=1}^{m} R_i$ is also a coherent
ring (by \cite[Theorem 2.4.3]{Glaz}).\\
Using the results of (strongly) Gorenstein flat modules
(\cite[Theorem 3.5]{BM}, \cite[Proposition 4.6]{BM2}, and
\cite[Theorem 3.7]{HH}),
 using  the fact that, if $I$ is an injective $R_1\times R_2$-module,
 then  $I=\Hom_{R_1\times R_2}(R_1\times R_2, I)
 =\Hom_{R_1\times R_2}(R_1, I)\times\Hom_{R_1\times R_2}(R_2, I)$
where each of $\Hom_{R_1\times R_2}(R_i, I)$, for $i=1,\ 2$,  is
an injective $R_i$-module, and by Lemma \ref{Prod-flat} below, the
argument is similar to the proof of Lemma
\ref{lem-produ-SGp}.\cqfd

\begin{lem}\label{Prod-flat}
Let $R_1\times R_2$ be a product of rings $R_1$ and $R_2$, and let
$F_{i}$ be an $R_i$-module for i=1, 2. Then, $\fd_{R_1\times
R_2}(F_1\times F_2)=\sup\{\fd_{R_1}(F_1), \fd_{R_2}(F_2) \}$.
\end{lem}
\proof We have
 $(F_1 \times F_{2})\otimes (R_1 \times 0)= (F_1 \times
F_{2})\otimes (R_1 \times R_2\, /\, 0\times R_2)= F_1\times 0.$
Then, since $R_1$ is a projective $R_1\times R_2$-module,
\cite[Chapter VI, Exercise 10, page 123]{H.Cartan} gives: $
\fd_{R_1}(F_1 \times 0)=\fd_{R_1}((F_1 \times F_{2})\otimes (R_1
\times 0)) \leq \fd_{R_1\times R_2}(F_1 \times F_{2}).$ The same:
$\fd_{R_2}(0\times F_2)\leq \fd_{R_1\times R_2}(F_1 \times
F_{2}).$\\
Thus, $\sup\{\fd _{R_1}(F_{1}), \fd
_{R_2}(F_{2})\}\leq \fd_{R_1 \times R_2}(F_1 \times F_{2}).$\\
Conversely, from  \cite[Chapter VI, Exercise 10, page
123]{H.Cartan}, we have:\\ $\fd_{R_1 \times R_2}(F_1 \times 0)\leq
\fd_{R_1 \times 0}(F_1 \times 0)= \fd_{R_1 \times 0}(F_1)$ and
$\fd_{R_1 \times R_2}(0 \times F_2)\leq \fd_{0 \times R_2}(0
\times F_2)= \fd_{0 \times R_2}( F_2).$ Therefore, $ \fd_{R_1
\times R_2}(F_1 \times F_{2})= \sup\{\fd _{R_1 \times
R_2}(F_{1}\times 0), \fd _{R_1 \times R_2}(0 \times F_{2})\} \leq
\sup\{\fd_{R_1}(F_{1}), \fd _{R_2}(F_{2})\}$, as
desired.\cqfd\bigskip

\noindent\textbf{Proof of Theorem \ref{GF-direct-product}.} Using
\cite[Theorem  4.11]{BM2}, \cite[Proposition 3.10]{HH}, and Lemma
\ref{lem-produ-SGf}, the argument is similar to the proof of
Theorem \ref{GP-direct-product}.\cqfd

Theorem \ref{GF-direct-product} allows us to construct a family of
non-Noetherian coherent rings $\{R_{i}\}_{i\geq 1}$ such that
$\Gwdim(R_i)=i$, $\Ggldim(R_i)>i$,  and $\wdim(R_i)=\infty$ for
all $i\geq 1$, as follows:

\begin{exmp}\label{Exm2}
Consider   a non-semisimple quasi-Frobenius  ring $R$ and a
semihereditary ring $S$ which is not hereditary.  Then, for every
positive
 integer $n$, we have:
 $\Gwdim(R\times S[X_1,X_2,...,X_n])=n+1,\quad \Ggldim(R\times S[X_1,X_2,...,X_n])>n+1,$ and
$
   \wdim(R\times S[X_1,X_2,...,X_n])=\infty.$
\end{exmp}
\proof Recall, at first, that every semihereditary ring $T$ is
stably coherent, that is the polynomial ring $T[X_1,...,X_j]$ is
coherent for every positive integer $j$ \cite[Corollary
7.3.4]{Glaz}.\\
The proof of the first and the last equality is similar to
the one of Example \ref{Exm1-non-Noeth}.\\
Now, assume,  by absurd,  that $\Ggldim(R\times
S[X_1,X_2,...,X_n])\leq n+1$.\\
 From Theorem \ref{GP-direct-product},
 $\Ggldim(S[X_1,X_2,...,X_n])\leq
n+1$. But, since $S$ is semihereditary, $\wdim(S)\leq 1$,   then
$\wdim(S[X_1,X_2,...,X_n])\leq n+1$. Thus, from \cite[Proposition
3.11]{BM2}, we have:
 $\gldim(S[X_1,X_2,...,X_n])=\Ggldim(S[X_1,X_2,...,X_n])\leq
 n+1.$
Hence, $\gldim(S)\leq 1$. But, this means that $S$ is hereditary
 which is absurd by hypothesis.\cqfd
\end{section}


\begin{thebibliography}{999}

\bibitem{BM}        D. Bennis and N. Mahdou; Strongly Gorenstein projective, injective, and flat modules, 
J. Pure Appl. Algebra (In Press). 
\bibitem{BM2}       D. Bennis and N. Mahdou; Global Gorenstein Dimensions, submitted for publication.  
                      Available from math.AC/0611358 v2 28 Feb 2007.
\bibitem{H.Cartan}  H. Cartan and S. Eilenberg; Homological Algebra, Princeton University Press, 1956.
\bibitem{LW}        L. W. Christensen; Gorenstein dimensions, Lecture Notes in Math., 1747, Springer, Berlin, (2000).
\bibitem{CFH}       L. W. Christensen, A. Frankild, and H. Holm; On Gorenstein projective, injective and flat dimensions - a functorial description with applications,  J. Algebra {\bf 302} (2006), 231--279.
\bibitem{Glaz}      S. Glaz; Commutative Coherent Rings, Springer-Verlag, Lecture Notes in Mathematics, 1371 (1989).
\bibitem{HH}        H. Holm; Gorenstein homological dimensions, J. Pure Appl. Algebra  {\bf 189} (2004), 167--193.
\bibitem{Iwa}       Y. Iwanaga; On rings with finite self-injective dimension, Comm. Algebra {\bf 7} (4) (1979), 393--414.
\bibitem{Kap}       I. Kaplansky; Commutative Rings, Allyn and Bacon, Boston, (1970).
\bibitem{Mah2001}   N. Mahdou;  On Costa's conjecture, Comm. Algebra {\bf 29}  (7) (2001) 2775-2785.
\bibitem{Rot}       J. J. Rotman; An introduction to homological algebra, Academic Press, New York, 1979.

    \end{thebibliography}
\end{document}